\documentclass[12pt]{amsart}
\usepackage{amsthm,amsmath,amsxtra,amscd,amssymb,xspace,esint}
\usepackage[cp1251]{inputenc}
\usepackage[T2A]{fontenc}

\numberwithin{equation}{section}

\def\beq{\begin{equation}}
\def\eeq{\end{equation}}
\def\p{\partial}
\def\G{\Gamma}
\def\g{\gamma}
\def\s{\sigma}

\def\bt{{\bf t}}
\def\e{\varepsilon}

\def\a{\alpha}
\def\b{\beta}
\def\l{\lambda}

\def\LL{{\mathcal L}}

\def\dim{{\rm dim}}
\def\res{{\rm res}}

\def\BA{Baker-Akhiezer\xspace}

\def\wt{\widetilde}
\def\wh{\widehat}

\def \matrix #1 {\left(\begin{array}{cc} #1 \end{array}\right)}
\newtheorem{thm}{Theorem}[section]

\newtheorem{cor}{Corollary}[section]
\newtheorem{lm}{Lemma}[section]
\newtheorem{rem}{Remark}

\begin{document}

\title[Curves with involution]{Characterizing Jacobians of algebraic curves with involution}

\author{Igor Krichever}
\address{Columbia University, New York, USA, and Skolkovo Institute for Science and Technology, and National Research University Higher School of Economics, Moscow, Russia}
\email{krichev@math.columbia.edu}

\begin{abstract}
We give two characterizations of Jacobians of curves with involution having fixed points in the framework of two particular cases of Welter's trisecant conjecture. The geometric form of each of these characterizations is the statement that such Jacobians are exactly those containing a shifted Abelian subvariety whose image under the Kummer map is orthogonal to an explicitly given vector.
\end{abstract}

\maketitle

\section{Introduction}

The problem of characterization of the Jacobians of curves among irreducible principally polarizaed Abelian varieties (ppav) is the famous Riemann-Schottky problem.
Over years of its more than 130 years history, quite a few geometrical characterizations of the Jacobians have been found (see surveys \cite{arb:expository,deb:expository,taimanov}).
None of them provides an explicit system of equations for the image of the Jacobian locus in the projective space under the level two theta imbedding.

The first effective solution of the Riemann-Schottky problem was obtained by T.Shiota \cite{shiota}, who proved the famous Novikov's conjecture:
{\it the Jacobians of smooth algebraic curves are precisely those indecomposable principally polarized abelian varieties (ppavs) whose theta-functions provide solutions to the Kadomtsev-Petviashvili (KP) equation}.

The strongest known characterization of a Jacobian variety in this direction, the so-called {\it Welters' (trisecant) conjecture\/}, formulated in \cite{wel1}: {\it an indecomposable ppav $X$ is the Jacobian of a curve if and only if its Kummer variety $K(X)$ has a trisecant line\/} was proved by the author in \cite{kr-schot,kr-tri}. The approach to its solution, proposed in \cite{kr-schot}, is general enough to be applicable to a variety of Riemann-Schottky-type problems. In \cite{kr-quad,kr-prym} it was used for a characterization of principally polarized Prym varieties. The latter problem is almost as old and famous as the Riemann-Schottky problem but is much harder.

\medskip
The goal of that work is to characterize the Jacobians of curves with involution. The curves with involution naturally appears as a part of algebraic-geometrical data defining solutions to integrable system with symmetries. Numerous examples of such systems  include the
Kadomtsev-Petviashvili hierarchies of type B and C (BKP and CKP hierarchies, respectively) introduced in  \cite{DJKM81,DJKM83}  and  the Novikov-Veselov hierarchy introduced in \cite{nv1,nv2}.

The existence of an involution of a curve is key in proving that the constructed solutions have the necessary symmetry.
The solutions corresponding to the same curve are usually parameterized by points of its Prym variety. In other words the existence of involution plus some extra constraints on the divisor of the Baker-Akhiezer function are {\it sufficient} conditions ensuring required symmetry. The problem of proving that these conditions are {\it necessary for two-dimensional} integrable hierarchies is much harder and that is a problem we address in that paper.

\begin{rem}
{\rm To the best of our knowledge from pure algebraic-geometri\-cal perspective the characterization problem of curves with involution in terms of their Jacobians has never been considered in its full generality. The only known to the author works in this directions are \cite{bv,poor,grush}. In a certain sense the setup we will consider -- the Jacobian and the Prym variety in it -- resembles the setup arising in the famous Schottky-Jung relations \cite{schot-jung}.}
\end{rem}

\medskip

Let $B$ be an indecomposable symmetric matrix with positive definite imaginary part.
It defines an indecomposable principally polarized abelian variety
$X=\mathbb C^g/\Lambda$, where  the lattice $\Lambda$ is generated by the basis vectors
$e_m\in \mathbb C^g$ and the column-vectors $B_m$ of $B$. Throughout the paper  $ \pi: \mathbb C^g \to X $ denotes the corresponding projection.

The Riemann theta-function $\theta(z)=\theta(z|B)$ corresponding to $B$
is given by the formula
\beq\label{teta1}
\theta(z)=\sum_{m\in \mathbb{Z}^g} e^{2\pi i(z,m)+\pi i(Bm,m)},\ \
(z,m)=m_1z_1+\cdots+m_gz_g
\eeq
The Kummer variety $K(X)$ is an image of the Kummer map
\beq\label{kum}
K:Z\in X\longmapsto
\{\Theta[\e_1,0](Z):\cdots:\Theta[\e_{2^g},0](Z)\}\in \mathbb{CP}^{2^g-1}
\eeq
where $\Theta[\e,0](z)=\theta[\e,0](2z|2B)$ are level two theta-functions
with half-integer characteristics $\e$.

A trisecant of the Kummer variety is a projective line which meets
$K(X)$ at least at three points. There are three particular cases of the characterization of
the Jacobians  by trisecants, corresponding to three possible configurations of
the intersection points $(a,b,c)$ of $K(X)$ and the trisecant:

(i) all three points coincide $(a=b=c)$,

(ii) two of them coincide $(a=b\neq c)$;

(iii) all three intersection points are distinct
$(a\neq b\neq c\neq a)$.

\medskip

Of course the first two cases can be regarded as degenerations of the general case~(iii).
However, when the existence of only one trisecant is assumed, all three cases are
independent and require its own approach. The approaches used in \cite{kr-schot,kr-tri} were based on the theories of three main soliton hierarchies (see details in \cite {kr-shiot}): the KP hierarchy for (i), the 2D Toda hierarchy for (ii) and the Bilinear Discrete Hirota Equations (BDHE) for (iii). Recently, pure algebraic proof of the first two cases of the trisecant conjecture were obtained in \cite{acp}.

\medskip
The main goal of this work is to give two characterizations of the Jacobians of curves with involution,
which distinguish such Jacobians within the framework of their characterizations given by cases (i) and (ii) above. Both of them are limited to the case of involutions having at least one fixed point, i.e. to two-sheeted {\it ramified} covers. The first, related to the KP theory, is limited by the obvious reason, since a curve with one marked point is used in constructing its solutions.

\begin{thm}\label{thm:main} An indecomposable principally polarized abelian variety $(X,\theta)$ is the Jacobian variety of a smooth algebraic curve $\G$ of genus $g$ with involution $\s:\G\to \G$  having at least one point fixed if and only if there exist $g$-dimensional vectors $U\neq\,0,V,A,\zeta $ and constants $\Omega_1,\Omega_2, b_1$ such that:

$(A)$ the equality
\beq
\left(\p_y-\p_x^2+u\right)\psi=0\,,\label{lax0}
\eeq
where
\beq\label{u0}
u=-2\p_x^2 \ln \theta (Ux+Vy+Z)+b_1,\ \ \ \ \psi=\frac{\theta(A+Ux+Vy+Z)}{\theta(Ux+Vy+Z)}\, e^{\Omega_1\,x+\Omega_2\,y},
\eeq
holds, for an arbitrary vector $Z$;

\medskip
\noindent
and

$(B)$ the intersection of the theta-divisor $\Theta=\{Z\in X\,\mid\, \theta(Z)=0\}$ with a shifted abelian subvariety $Y\subset X$ which is the Zariski closure of $\pi(Ux+\zeta) \subset X$ is {\rm reduced} and the equation
\beq\label{C}
\p_V\theta|_{\Theta\cap Y}=0
\eeq
holds.

Moreover, the locus $\Pi$ of points $\zeta\in X$ for which the equation \eqref{C} holds is the locus
of points for which the equation $\zeta+\sigma(\zeta)=2P+K\in X $, where $K$ is the canonical class,
holds.
\end{thm}
\begin{rem} {\rm Note that if $V\neq 0$ then \eqref{C} implies that $\dim\, Y<g$, since for any nonzero vector $V$ the restriction $\p_V\theta|_\Theta\neq 0$. The case $V=0$ when \eqref{C} is automatically satisfied is the case of hyperelliptic curves.}
\end{rem}
The condition $(A)$ is one of three equivalent forms of the characterization of the Jacobains among ppav proved in \cite{kr-schot}. The direct substitution of the expression (\ref{u0}) into equation (\ref{lax0}) and the use of the addition formula for the Riemann theta-functions
\beq\label{ThetaQuad}
\theta(z+w)\theta(z-w)=\sum_{\e\in((1/2)\mathbb Z/\mathbb Z)^g}\Theta[\e,0](z)\Theta[\e,0](w)\,.
\eeq
shows the equivalence of $(A)$ to the condition that for all theta characteristics $\e\in (\frac 12\mathbb Z/\mathbb Z)^g$ the equations
\beq\label{gr0}
\left(\p_V-\p_U^2-2\Omega_1\,\p_U+(\Omega_2-\Omega_1^2)\right)\, \Theta[\e,0](A/2)=0
\eeq
hold. (Here and below $\p_U$, $\p_V$ are the derivatives along the vectors $U$ and $V$, respectively). Equations \eqref{gr0} means that the image of the point $A/2$ under the Kummer map is an inflection point (case~(i) of Welters' conjecture).

There are two other equivalent forms of the condition $(B)$, which, in particular, give its {\it geometric form}.  The first one is:

\medskip
$(C)$ {\it there is a vector $W$ and a constants $\Omega_3,b_3$ such that the equality
\beq\label{a}
\left(\p_t-\p_x^3+\frac 3 2 u\p_x+\frac 34 u_x+b_3\right)\psi=0\,,
\eeq
where
\beq\label{a0}
u=-2\p_x^2 \ln \theta (Ux+Wt+\zeta)+b_1,\ \ \ \ \psi=\frac{\theta(A+Ux+Wt+\zeta)}{\theta(Ux+Wt+\zeta)}\, e^{\Omega_1\,x+\Omega_3\,t},
\eeq
holds.}

The fact that the condition $(C)$ holds for curves with involution is known.  In the CKP theory \cite{DJKM81,DJKM83} equation \eqref{a} plays the same role as equation \eqref{lax0} in the KP theory. Namely, both equations define the first flows of the corresponding hierarchies (see details in \cite{kz}).

In \cite{kn,kz} equation \eqref{a} was obtained and used in another but equivalent form:

\medskip
$(C')$ {\it there is a constant $b_2$ such that the equality
\beq\label{b1}
\p_U\p_V\ln \theta|_{\widehat Y}=b_2
\eeq
 holds on $\widehat Y$ which is Zariski closure of $\pi(Ux+Wt+\zeta)\subset X$}.

\begin{rem}\label{rm:bs} {\rm As will be seen in what follows, if there is a set of vectors and constants for which the above conditions are satisfied, then there is a family of such sets. One of the consequences of this is that without loss of generality, any (but only one) of the constants $b_i$ above can be set equal to $0$}.
\end{rem}

\medskip

The addition theorem \eqref{ThetaQuad}  implies that \eqref{b1} is equivalent to the condition that the vector $\left(\p_U\p_V K(0)-b_2K(0)\right)$ is {\it orthogonal to the image under the Kummer map $K(\Pi)$ of the shifted abelian subvariety} $\widehat Y$:
\beq\label{ort}
\sum_{\e\in((1/2)\mathbb Z/\mathbb Z)^g}\left(\p_U\p_V\Theta[\e,0](0)-b_2\Theta[\e,0](0)\right)\Theta[\e,0](z)=0,\ \ z\in \widehat Y
\eeq
whence follows the condition {\it of a kind of flatness} of the image under the Kummer map of the shifted Prym subsubvariety $\Pi \subset X $, that is, $ K (\Pi) $ lies in a proper (projective) linear subspace.

\medskip
Equations \eqref{lax0} and \eqref{b1} are used for the proof of $(B)$, which is the strongest form of the characterization the Jacobians of curves with involution.  The implication $(B)\Rightarrow (C)$  is not  by all means direct and comes only as a result of the proof the theorem. Namely, $(B)$ is what we {\it really}
use in the proof that the corresponding curve is a curve with involution. The latter, as  it was mentioned above implies $(C)$.

\medskip
The explicit meaning $(B)$ is as follows. As shown in \cite{shiota,Fay} the affine line $Ux+Z$ is not contained in $\Theta$ for any vector $Z$. Hence, the function $\tau(x,y):=\theta(Ux+Vy+z), z\in Y$ is a {\it nontrivial} entire function of $x$. The assumption that $\Theta\cap Y$ is reduced means that zeros $q(y)$ of $\tau$ considered as a function of $x$ (depending on $y$) are generically simple, $\tau(q(y),y)=0,\, \tau_x(q(y),y)\neq 0$. Then \eqref{b1} is the equation
$\p_y\, q|_{y=0}=0$.

\medskip

The second characterization of the Jacobians of curves with involution is related to the 2D Toda theory. A priory, unlike the KP case, there is no obvious reason why it is not applicable to all types of involution including unramified covers. It turned out that there is an obstacle for the case unramified covers and our second theorem also gives a characterization of the Jacobians of curves with involution {\it with} fixed points.

\begin{thm}\label{main2}
An indecomposable, principally polarized abelian variety $(X,\theta)$
is the Jacobian of a smooth curve of genus $g$ with involution having fixed points if and only if
there exist non-zero $g$-dimensional vectors
$U\neq A \, (\bmod\,  \Lambda),\, V, \zeta$, constants $\Omega_0,\Omega_1,b_1$ such that the following two conditions are satisfied.

$(A)$  The differential-functional equation
\beq\label{laxd}
\left(\p_y-T-u\right)\psi=0, \ \ T=e^{\p_x}
\eeq
where
\beq\label{u21}
u=b_1+(T-1)\p_y\ln\theta (xU+yV+Z)
\eeq
and
\beq\label{p21}
\psi=\frac{\theta(A+xU+yV+Z)}{\theta(xU+yV+Z)}\, e^{x\Omega_0+y\Omega_1}
\eeq
holds for an arbitrary vector $Z$.

\medskip
$(B)$ $(i)$ The intersection of the theta-divisor with the shifted Abelian variety $Y$, which is a closure of $\pi(Ux+\zeta)$, is reduced and is not invariant under the shift by $U$, $\Theta\cap Y\neq (\Theta+U)\cap Y$, and $(ii) $ the equation
\beq\label{Cd}
\left((\p_V\theta(z))^2+b_2\,\theta(z+U)\theta(z-U)\right)|_{z\in\Theta\cap Y}=0,
\eeq
where $b_2\neq 0$ is a constant, holds.

Moreover, the locus of the points $\zeta\in X$ for which the equation \eqref{Cd} holds is the locus of point for which
the equation $\zeta+\zeta^\s=K+P_1+P_2\in J(\G)$, where $(P_1, P_2)$ are points of the curves permuted by $\s$ and such that $U=A(P_2)-A(P_1)$, is satisfied.

\end{thm}

\begin{rem} Under the assumption that $U$ spans an elliptic curve Theorem \ref{main2} was proved in \cite{kz2}
\end{rem}
The condition $(A)$ is one of three equivalent forms of the characterization of the Jacobians proved in \cite{kr-tri}.
It is equivalent to the condition (case $(ii)$ of the Welter's conjecture):

{\it the equations
\beq\label{gr1}
\p_{V}\Theta[\e,0]\left((A-U)/2\right)-e^{p}\Theta[\e,0]\left((A+U)/2\right)
+E\Theta[\e,0]\left((A-U)/2\right)=0,
\eeq
are satisfied for all $\e\in  Z_2^g$}.

\medskip

In the course of proving that condition $(B)$ holds for ramified double covers, we first prove that:

{\it there is a vector $W$ and constants $\Omega_2, b_2$ such that the differential-functional equation

\beq\label{semi1}
(\p_t-T-w_1-w T^{-1})\psi=0
\eeq
where
\beq\label{w1}
w_1=b_1+\frac 12 (T-1)\p_t\ln {\theta}, \quad w=b_2\, \frac{T\theta T^{-1}\theta}{\theta^2}
\eeq
$$\theta=\theta (xU+tW+\zeta)$$
and
\beq\label{psid}
\psi=\frac{\theta(A+xU+tW+\zeta)}{\theta(xU+tW+\zeta)}\, e^{x\Omega_0+t\Omega_2},
\eeq
holds}.

\medskip
The equivalent form of the statement above is the statement that for the ramified double covers

\medskip
{\it there is a constant $b_3$ such that the equality
\beq\label{bd1}
\theta^2(z)\p^2_V\ln\theta(z)-b_2\, \theta(z+U)\theta(z-U)|_{z\in\widehat Y}=b_3\,,
\eeq
where $\widehat Y$ which is  Zariski closure of $\pi(Ux+Wt+\zeta)$, holds. }

\medskip

The addition theorem \eqref{ThetaQuad}  implies that \eqref{bd1} is equivalent to the condition that the vector $(2\p_V^2 K(0)-b_2K(U)-b_3K(0))$ is orthogonal to the image under the Kummer map of the abelian subvariety $\widehat Y$, i.e. the equation
\beq\label{ortd}
\sum_{\e\in((1/2)\mathbb Z/\mathbb Z)^g}(2\p^2_V\Theta[\e,0](0)-b_2\,\Theta[\e,0](U)-b_3\Theta[\e,0](0)))\Theta[\e,0](z)=0
\eeq
with $z\in \widehat Y$, holds.

\smallskip
Equations \eqref{semi1} and \eqref{bd1} are analogs of the conditions $(C)$ and $(C')$ in the flex case. But unlike the latter, we do not claim that they are equivalent to $(B)$, because we came short in proving that they imply $(B)$, and will use additional arguments for the last step in the proof.  The condition $(B)$ is what we really use for the proof of "if" part of the theorem.

\begin{rem} \label{rm:W} {\rm As we shall see below, the case $2U\in \Lambda$ and $W=0$ in \eqref{w1} and \eqref{psid} corresponds to hyperelliptic curves, which are curves with involution}.
\end{rem}

We conclude the introduction by saying that it is tempting to see if the third characterization associated with the fully discrete BDHE hierarchy might be applicable to the case of unramified covers.

\smallskip
{\bf Aknowledments}. The author would like to thank Enrico Arbarello, Sam Grushevsky and Anton Zabrodin for very useful and inspiring discussions.

\section{Preliminaries}

\subsection*{The KP and CKP hierarchies}

The Kadomtsev-Petviashvily (KP) hierarchy is one of the most fundamental in the modern
theory of integrable systems. It has at least three well-known
de\-fi\-ni\-ti\-ons/rep\-re\-sen\-ta\-ti\-ons usually called: the Zakharov-Shabat form \cite{ZS75}, the Sato form
\cite{sato} and the Hirota bilinear equation form \cite{DJKM83,JimboMiwa}.

In the Sato form it is defined as a system of commuting flows
on the space of sequences $(u_1(x),u_2(x),\ldots) $
of functions of one variable $x$, which can be identified
with the space of pseudo-differential operators of the form
\beq\label{satokp}
\LL=\p_x +u_1\p_x^{-1}+u_2 \p_x^{-2}+\ldots
\eeq
The flows are defined by the Lax equations
\beq\label{kp3}
\p_{t_k} \LL=[B_k, \, \LL], \quad B_k:=\LL^k_+ \quad k=1,2,3, \ldots
\eeq
where $(\cdot)_+$ stands for the differential part of a pseudo-differential operator. The commutativity of flows implies that the equations
\beq\label{kp10}
\p_{t_l}B_k-\p_{t_k}B_l+[B_k, B_l]=0
\eeq
are satisfied for all pairs $k,l$.  For $k=2,l=3$  operator equation \eqref{kp10} where
\begin{eqnarray}\label{kp149}
B_2&=\p_{x}^2-u,\ \ \ \ \quad \quad\quad  &u=-2u_1\\
B_3&=\p_{x}^3-\frac 32 u\p_{x}-w, \quad &w= -\frac 32 u_x-3u_2 \label{kp148}
\end{eqnarray}
after change of notation $t_1=x, t_2=y,\, t_3=t$, is equivalent to the system of two equations
\beq\label{kp100}
4w_x=3u_{xx}+3u_y, \quad 4w_y=(4u_t+6uu_x-u_{xxx})_x+3u_{xy}
\eeq
Eliminating $w$ from the system one gets the original KP equation for the remaining function $u$
\beq\label{kporig}
3u_{yy}=(4u_t+6uu_x-u_{xxx})_x
\eeq

In \cite{DJKM81} an infinite integrable hierarchy of partial differential was introduced and  called the
Kadomtsev-Petviashvili hierarchy of type C (CKP). It is a hierarchy of commuting flows that are the restriction of the flows of the KP hierarchy corresponding
to ``odd'' times ${\bf t}_{\rm o}=\{t_1, t_3, t_5, \ldots \}$
onto the space of {\it anti self-adjoint} pseudo-differential operators $\LL$ of the form (\ref{satokp}), i.e. such that
\beq\label{ckp01}
\LL^*=-\LL,
\eeq
where $()^*$ means the formal adjoint defined by the rule $\Bigl (f(x)\circ \p_x^{m}\Bigr )^{*}=(-\p_x)^m \circ f(x)$.

\begin{rem}
{\rm Note, that if \eqref{ckp01} is satisfied then the operator $B_3$ in \eqref{kp148} is of the form  \eqref{a}}.
\end{rem}
\medskip

In \cite{kz} the CKP hierarchy was characterized in terms of the
KP tau-function which is a function of KP "times", $\tau(t_1,t_2,\ldots)$. More precisely, each solution of the CKP hierarchy has a unique extension to the solution of the full KP hierarchy via the flows (\ref{kp3}) with even $k$ (which obviously do not preserve constraint (\ref{ckp01})). Such solution is naturally to call {\it KP extension of the solution to the CKP hierarchy}. In \cite{kz} it was proved that the KP tau-function is the tau-function of such a solution if and only if the equation
\beq\label{int3}
\p_x\p_{y}\log \tau\Bigl |_{{\bf t}_{\rm e}=0}=0
\eeq
holds for all ${\bf t}_{\rm o}$, where all ``even'' times ${\bf t}_{\rm e}=(t_2=y,t_4,\ldots)$ are set equal to zero.

\medskip
Theorem \ref{thm:main} is a stronger version of that result in the algebraic-geometric setting when $\tau$ function coincides with the theta-function up to a factor which is the exponent of a quadratic in times form. Namely it states, that it is enough to require that equation \eqref{int3} holds only for the first two times $t_1=x, t_3=t$ of the CKP hierarchy (compare with \eqref{b1}).
\subsection*{The Baker-Akhiezer function.} Recall that a smooth genus $g$ algebraic curve $\G$ with fixed local coordinate $k^{-1}(p)$  in the neighborhood of a point $P\in \Gamma$, $k^{-1}(P)=0$ and a generic effective degree $g$ divisor $D={\gamma}_{1}+{\cdots}+{\gamma}_{g}$ defines the Baker-Akhiezer function $\psi(\bt ,p)$ which is a function of complex variables $\bt=(t_1,t_2,t_3,\ldots)$ (it is always assumed that only finite number of $t_i\neq 0)$ and $p\in \Gamma$. For fixed $\bt$ it is defined as a unique function of $p\in \Gamma$ with the following analytic properties:

$1^0$. Outside  $P$ the singularities of $\psi$ are poles at points of the divisor $D$ of order not greater
then the multiplicity of the point in $D$, i.e., $(\psi)+D\geq 0$.

$2^0$. In the neighborhood of $P$ the function $\psi$ has the form
\beq
\psi(\bt,k)=e^{\sum_i k^i\, t_i}
\left(1+\sum_{s=1}^{\infty}\xi_{s}(\bt)\, k^{-s}\right)\, , \  k=k(p), \label{psi}
\eeq
In order to present an explicit formula for $\psi$ in terms of the Riemann theta function we first choose a symplectic basis of $a$-  and $b$-cycles on $\G$. Then define a basis of normalized holomorphic differentials $ \omega_k , \, k = 1,\ldots, g,$ and the matrix $B$ of
their $b$-periods
$$\oint_{a_k}\omega_\ell=\delta_{k\,\ell},\quad B_{k\, \ell}=\oint_{b_\ell}\omega_k=B_{\ell\, k}\,,$$
and the Riemann-theta function by formula \eqref{teta1}

\medskip
Denote by $A(p)$ the vector (depending on a path of integration) with coordinates  $A_k(p)=\int_{P}^p \omega_k$ and by $\Omega_i(p), i=1,2,3,\ldots,$ the Abelian integral $\Omega_{i}=\int_{P}^p d\Omega_{i}$ where $d\Omega_i$ is the normalized (i.e. having zero $a$-periods) meromorphic differential whose the only singularity is at $P$ of the form $d\Omega_i =dk^i(1+O(k^{-i-2}))$.

The definition of $\Omega_i$ needs a clarification since $d\Omega_i$ has the pole at $P$. In the definition of $\Omega_i$ it is assumed that its branch in the neighborhood of $P$ is fixed such that there is no constant term in the expansion
\beq\label{constb}
\Omega_i=k^{i}+\sum_{s=1}^\infty a^{(i)}_s k^{-1}+O(k^{-2})
\eeq
and then extended analytically along the path. It is assumed that the paths in the definition of $A(p)$ and $\Omega_i(p)$ are the same.

\begin{lm}[\cite{kr1,kr2}]\label{lm:BA} (i) The Baker-Akhiezer function defined above equals
\beq\label{formula}
\psi(\bt,p)=\frac{\theta({A}(p)+\sum_i t_iU_i+Z)\,\theta({Z})}
{\theta(\sum_i t_iU_i+{Z})\, \theta({A}(p)+Z)} \ e^{\sum_i t_i\Omega_{i}(p)}\ ,
\eeq
where  $U_i$ are vectors with the coordinates
\beq\label{U}
U_i^k=\frac{1}{2\pi i}\oint_{b_k}d\Omega_i,
\eeq
and
\beq\label{Z}
Z=-\sum_s A(\g_s)+{\mathcal K},
\eeq
where ${\mathcal K}$ is the vector of Riemann constants.

(ii) The BA function $\psi$ satisfies the equations
\beq\label{kp150}
\left(\p_{t_k}-B_k\right)\psi=0, \quad k=1,2,\ldots
\eeq
where $B_k$ is a monic differential operator in $x$ of order $k$.
\end{lm}
\begin{rem} {\rm The compatibility conditions of equations \eqref{kp150} is the KP hierarchy in the Zakharov-Shabat form. From the definition of $B_k$ in \cite{kr1} it is easy to show that $B_k=\LL^k_+$ where $\LL$ is a unique pseudo-differential operator such that the equation
\beq\label{kp151}
\LL\psi=k\psi
\eeq
holds. The compatibility conditions of \eqref{kp150} and \eqref{kp151} are the Lax equations \eqref{kp3}.}
\end{rem}

The substitution of (\ref{psi}) into \eqref{kp150} with $B_2$ of the form \eqref{kp149}  gives the equations
\beq\label{eqxi}
\p_y\xi_s-2\p_x\xi_{s+1}-\p_x^2 \xi_s+u\xi_s=0,\ \ s=0,\ldots
\eeq
The first of them with $s=0$ gives an expression of $u$ in terms of $\xi_1$. Namely, $u=2\p_1\xi_1$.

Similarly, the substitution \eqref{psi} into \eqref{kp150} with $B_3$ of the form \eqref{kp148}
gives the equations:
\beq\label{eqxi1}
\p_t\xi_s-3\p_x\xi_{s+2}-3\p^2_x\xi_{s+1}-\p_x^3 \xi_s+\frac 32 u\xi_{s+1}+\frac 32 u\p_x \xi_{s}+w\xi_s=0.
\eeq
The first of them with $s=0$ gives the expression for $w$.

\medskip
The explicit theta-functional formulas for $u$ and $w$ are obtained by expansion of \eqref{formula} near $P$. From the bilinear Riemann identities for periods of Abelian differentials it follows that the expansion of the Abel map $A(p)$ near $P$ has the form
\beq\label {Ak}
A(k)=-\sum_{i=1} \frac 1i \,U_i k^{-i}, \ k=k(p)
\eeq
where the vectors $U_i$ are given by formula \eqref{U}. From \eqref{formula} and \eqref{u29} it follows that
\beq\label{u29}
u=2\p_x\xi_1=-2\p_x^2 \ln \theta\left(\sum_it_i U_i+Z\right)+2a^{(1)}_1
\eeq
where the constant $a^{(1)}_1 $ is defined in \eqref{constb}.

Computing the next coefficient of expansion of \eqref{formula} one gets
\beq\label{w}
w=\frac 34 u_x-\frac 32 \p_x\p_y\ln \theta\left(\sum_it_i U_i+Z\right)+b_3
\eeq

Although for completeness we defined above the Baker-Akhiezer function, depending on the full set of times of the KP hierarchy, in the future we will restrict ourselves to its dependence only on the first three times, setting $t_i = 0, i> 3$. Recall that above we have already denoted these times by $x =t_1, y=t_2, t =t_3 $. Below we will use notation $ U, V, W $ for the vectors $ U_1, U_2, U_3 $, respectively.

\subsubsection*{Bloch properties of the Baker-Akhiezer function}

For further use we elaborate more on monodromy properties of the Baker-Akhiezer function along the subvariety $Y\subset J(\G)$. Recall that $Y$ is Zariski closure of $\pi(xU+\zeta)$ where $\pi:\mathbb C^g  \longmapsto J(\G)=\mathbb C^g/\Lambda$ is the projection.

The preimage $\widetilde Y=\pi^{-1}(Y)$ is an affine subspace.  A choice of a vector
$\zeta\in \widetilde Y$ identifies it with the image of an imbedding map
\beq\label{iota}
\iota: \mathbb C^d \hookrightarrow \mathbb C^g
\eeq
Let $\Lambda_U\subset \Lambda$ be the sublattice of the Jacobian lattice such that under this identification
\beq\label{lambda}
Y= \mathbb C^d/\Lambda_U
\eeq
For further use define the vector $\eta\in \mathbb C^d$ by the equation
\beq\label{ueta}
\iota(\eta)=U
\eeq

\medskip
Consider the function
\beq\label{phi}
\phi(z,k):= \frac{\theta({A}(k)+\iota(z))\,\theta(\zeta)}
{\theta(\iota(z))\, \theta(A(k)+\zeta)}\, e^{(\ell(k),z)}\,, \quad z\in \mathbb C^d,
\eeq
where $A(k)$ is defined in \eqref{Ak} and $\ell(k)$ is a formal series
\beq\label{ell}
\ell(k)=\sum_{s=1}^\infty \ell_{s}k^{-s}, \quad \ell_{s}\in \mathbb C^d,
\eeq
such that the equation
\beq\label{oell}
(\ell(k),\eta)=\Omega_1(k)-k\,,
\eeq
holds.

\smallskip
A simple comparison of \eqref{formula} and \eqref{phi} shows that the coefficients $\xi_s(x,0)$ of the Baker-Akhiezer function $\psi(x,0,k)$ expansion at $P$ are of the form
\beq\label{psiphi1}
\xi_s(x,0)=\varphi_s(x\eta)
\eeq
where $\varphi_s(z)$ are coefficients of the expansion of $\phi$,
\beq\label{phiexp}
\phi(z,k)=1+\sum_{s=1}^\infty \varphi_s(z) k^{-s}
\eeq
The latter are of the form
\beq\label{phi11}
\varphi_s(z)=\frac {\tau_s(z)}{\tau(z)} , \quad \tau(z):=\theta(\iota(z))
\eeq
where $\tau_s(z)$ is a holomorphic function of $z$.

The monodromy properties of the theta-function imply that under the shift of the arguments by a vector $\lambda\in \Lambda_U$  the function $\phi$ gets transformed to
\beq\label{phitrans}
\phi(z+\lambda,k)=\phi(z,k)\rho_\lambda(k), \ \ \quad \rho_\lambda(k)=1+\sum_{s=1}^\infty\rho_{\lambda,s}k^{-s}
\eeq
Equation \eqref{phitrans} is equivalent to the equations
\beq\label{phitrans1}
\varphi_s(z+\lambda)=\sum_{j=0}^s \varphi_{s-j}(z)\rho_{\lambda, j}
\eeq

\begin{lm}\label{lm:periodicity} Let $\l_1,\ldots,\l_d$ be a set of linear independent vectors in $\Lambda_U$. Then there is a unique coordinate $k^{-1}(p)$ in the neighborhood of $P$ and a unique  linear form $\ell(k)$ satisfying \eqref{oell} such that $\phi(z,k)$ is invariant under the shift by these vectors, i.e.,
\beq\label{periodicity}
\phi(z+\lambda_i,k)=\phi(z,k) \quad \Leftrightarrow \quad \rho_{\lambda_i}(k)=1, \ \ i=1,\ldots, d.
\eeq
\end{lm}
For the proof it is enough to note, that the differential $d\Omega_1$ depends only on the first jet
of the local coordinate, i.e., it does not change under the change of the coordinate of the form $k'=k+O(k^{-1})$, but the coefficients of his expansion \eqref{constb} {\it do depend} on a choice of the local coordinate. Conversely, for  any given set of coefficients $a^{(1)}_{\,s}$ equation \eqref{constb} with $i=1$ can be regarded as the definition of the corresponding formal local coordinate.

The vectors $\ell_s$ and constants $a^{(1)}_{\,s}$ are defined recurrently. On each step we defined first $\ell_s$ by a system of nonhomogeneous linear equations which are needed for the equation $\rho_{\lambda_i,s}=0$ to be satisfied, and then define $a^{(1)}_{\,s}$  by the equation $a^{(1)}_{\,s} =(\ell_s,\eta)$.

\subsection*{The dual Baker-Akhiezer function}
A notion of duality is defined first for a generic effective degree $g$ divisor $D=\g_1+\ldots+\g_g$. For any such divisor there is a unique meromorphic differential $d\Omega$ having second order pole at the marked point $P$,
\beq\label{coeff}
d\Omega=dk\left(1+\sum_{s=2}^\infty c_s k^{-s}\right)
\eeq
vanishing  at the points of the divisor $D$ with the multiplicity at least equal to the multiplicity of the point in $D$. The zero divisor of $d\Omega$ is of degree $2g$. Let $D^*$ be the effective degree $g$ complimentary divisor, i.e. $D+D^*=K+2P\in J(\G)$, where $K$ is the canonical class.

\medskip
The dual BA function is defined as a unique function $\psi^*(x,y,p)$ with the following analytic properties with respect to $p\in \Gamma$:

$1^0$. Outside  $P$ the singularities of $\psi^*$ are poles at the divisor $D^*$, $(\psi*)+D^*\geq 0$;

$2^0$. In the neighborhood of $P$ the function $\psi$ has the form
\beq
\psi^*(x,y,k)=e^{-k\,x-k^2\,y}
\left(1+\sum_{s=1}^{\infty}\xi^*_{s}(x,y)\, k^{-s}\right)\, , \  k=k(p), \label{psidual}
\eeq
The explicit formula for $\psi^*$ in terms of the Riemann theta-function is:
\beq\label{formuladual}
\psi^*(x,y,k)=\frac{\theta({A}(p)-xU-yV-Z)\,\theta({Z})}
{\theta(xU+yV+Z)\, \theta({A}(p)-Z)} \ e^{-x\,\Omega_{2}(p)-y\, \Omega_{3}(p)}\ ,
\eeq

From the definition of the dual BA function it is easy to see that the equations
\beq\label{nomega}
\res_P\left(\p_x^i\psi^*(x,0,p)\p_x^m \psi(x,0,p)) \, \right) d\Omega=0, \quad m,i=0,1,\ldots
\eeq
hold. Indeed, the differential in the left hand side of the equation is holomorphic away of $P$ since the poles of $\psi$ and $\psi^*$ cancel with zeros of $d\Omega$. The essential singularities of $\psi$ and $\psi^*$ at $P$ cancel each other. Hence, the differential in l.h.s. is a meromorphic differential with the only pole at $P$.

Computing the residue with $i=0$ in terms of the coefficients of the expansions \eqref{psi} and \eqref{psidual} we get the system of equations
\beq\label{reseq}
\xi_{m+1}+\xi_{m+1}^*+c_{m+1}=h_{m+1}(\xi_1,\xi_1^*,\ldots,\xi_m,\xi_m^*; c_2,\ldots,c_m)
\eeq
where $h_{n+1}$ are some explicit differential polynomials in $\xi_s$, linear in $\xi_s^*, \, s\leq n,$ and in the coefficients of the  expansion \eqref{coeff} of $d\Omega$.
Equations \eqref{reseq} recurrently express coefficients of the dual function expansion in terms of the coefficients $\xi_s$ and $c_s$.

\medskip
The dual BA function satisfies the equation that is formally adjoint to \eqref{lax0}:
\beq\label{laxdual}
(-\p_y-\p_x^2+u(x,y))\psi^*(x,y,p)=0
\eeq
To prove \eqref{laxdual} it is enough to note that the same arguments as in the proof of \eqref{lax0} show that $\psi^*$ satisfies \eqref{laxdual}
with potential equal to  $v=-2\p_x\xi_1^*$. Equation \eqref{reseq} for $m=1$, i.e. the equation
\beq\label{res1}
\xi_1+\xi^*_1=0
\eeq
implies $v=u$ and \eqref{laxdual} is proved.

The substitution of \eqref{psidual} into \eqref{laxdual} gives a system of equations
\beq\label{eqxidual}
-\p_y\xi^*_s+2\p_x\xi^*_{s+1}-\p_x^2 \xi^*_s+u\xi^*_s=0,\ \ s=0,\ldots
\eeq

\medskip

Introduce the function dual to $\phi$:
\beq\label{phidual}
\phi^*(z,k):= \frac{\theta({A}(k)-\iota(z))\,\theta({\zeta})}
{\theta(\iota(z))\, \theta(A(k)-\zeta)}\,e^{-(\ell(k),z)},
\eeq
and the coefficients $\varphi_s^*(z)$ of its expansion. From \eqref{formuladual} it follows that
\beq\label{psiphi1star}
\xi^*_s(x,0)=\varphi^*_s(x\eta)
\eeq
The functions $\varphi_s^*$ have the form
\beq\label{phidual11}
\varphi^*_s(z)=\frac {\tau^*_s(z
)}{\tau(z)},
\eeq
where $\tau_s^*$ is a holomorphic function of $z\in \mathbb C^g$.
They  satisfy the following monodromy properties
\beq\label{phitrans2}
\varphi^*_s(z+\lambda)=\sum_{j=0}^s \varphi^*_{s-j}(z)\rho^*_{\lambda, j}
\eeq
where $\rho^*_{\lambda, s}$ are the coefficients of the expansion of $\rho_\lambda^*(k)=\rho_\lambda^{-1}(k)$.

\section{Proof of Theorem \ref{thm:main}}

To begin with let us show that if a curve $\Gamma$ admits an involution under which $P$ is fixed and $k$ is odd,
$k(p)=-k(\sigma(p))$, then the conditions $(C)$ and $(C')$ of the theorem are satisfied for $\zeta\in J(\G)$ such that
\beq\label{zeta}
\zeta+\zeta^\sigma=K+2P\in J(\G)
\eeq
In other words for $\zeta$ that is the image under the Abel transform of a divisor $D$ such that $D+\sigma(D)$ is the zero-divisor of a meromorphic differential $d\wt\Omega$ with the only pole of second order at $P$.

Consider the differential $(\psi(x,0,t,\sigma(p))\p_x\psi(x,0,t, p))d\wt\Omega$. It is a meromorphic differential on $\Gamma$ with the only pole at $P$. Hence, it has no residue at $P$. Computing this residue in terms of the coefficients of the expansion (\ref{psi}) we get
\beq\label{resPP}
2\xi_2(x,0,t)-\xi_1^2(x,0,t)+\p_x\xi_1(x,0,t)+c_2=0
\eeq
where $c_2$ is a constant defined by the Laurent expansion of $d\wt\Omega$ at $P$.

Taking the $x$-derivative of (\ref{resPP}) and using (\ref{eqxi}) with $s=1$ we get the equation
\beq\label{respp}
0=\p_x(2\xi_2(x,0,t)-\xi^2_1(x,0,t)+\p_x\xi_1(x,0,t))=\p_y\xi_1|_{y=0}
\eeq
Recall, that $\xi_1=-\p_x\ln\theta+a^{(1)}_{\,1}x+a^{(2)}_{\,1} y$. Therefore, equation \eqref{b1} with $b_2=a^{(2)}_{\,1}$ is satisfied for all $\zeta$ satisfying \eqref{zeta}. Then from \eqref{w} it follows that $w(x,0,t)=\frac 34 u_x+d$, i.e equation \eqref{a} holds.

The formulae for $\psi$ in \eqref{u0} and \eqref{a0} are just the evaluation of the formula \eqref{formula} for the Baker Akhiezer function at any $p\in \G$.

\medskip
Our next goal is to prove $(C)\Rightarrow (B)$. As shown in \cite{shiota,Fay} the theta divisor $\Theta$ in the Jacobian $J(\G)$ of an a smooth algebraic curve $\G$ does not contain a $\p_U$ invariant locus. In other words: for any vector $Z\in J(\G)$ the function $\tau(x;y,t):=\theta(Ux+Vy+Wt+Z)$ is a {\it non-zero} entire function of the variable $x$ depending on the variables $y$ and $t$ and $Z$. For brevity we omit explicit indication on $Z$-dependence. The zeros $q$ of $\tau$ in $x$ (depending on $y,t,Z$) correspond to intersection points of $Ux+Vy+Wt+Z$ with the theta-divisor $\Theta$. Since $\Theta$ is reduced and $U$-direction is transverse to $\Theta$ on an Zarisiki open set, generic zeros of $\tau$ are simple. It might be not the case when we consider $z\in \widehat Y$.

The multiplicity of a zero of an entire function in $x$ depending on parameters  is an upper-continuous function of the parameters. Hence, there is a constant $N$ such that for $z$ in an open everywhere dense set $(\Theta\cap Y) ^0\subset \Theta\cap \widehat Y$ the multiplicity of the zero $q(t)$ in the neighborhood of $x=0$ of the function $\theta(Ux+Wt+z), z\in (\Theta\cap Y) ^0$ for sufficiently small $t$ is exactly $N$. Then, by the implicit function theorem $q(t)$ is a smooth function of $t$. Then the Laurent expansion of $u=-2\p_x^2 \ln\tau$ is of the form
\beq\label{uL} u(x,0,t)=\frac {2N}{(x-q(t))^2}+v(t)+w(t)(x-q(t)+\ldots
\eeq
whose coefficients $v,w,\ldots$  are smooth functions of $t$.

The function $\psi(x,0,t)$ given by formula \eqref{u0} is a meromorphic function of $x$ with a pole at $x=q(t)$ of order $n\leq N$ (the numerator might have a zero at $q(t)$), i.e. its Laurent expansion of $\psi$ at $q$ has the form
\beq\label{psiL}
\psi(x,0,t)=\frac {\alpha(t)}{(x-q(t))^n}+O\left((x-q(t))^{-n+1}\right), \quad \alpha\neq 0
\eeq
Substitution of \eqref{psiL} into equation \eqref{a} and calculation of the coefficient at the leading term of the Laurent expansion of the left-hand side, which is of order $-n-3$ gives the equation
\beq\label{N1}
3N=n(n+2)
\eeq
Now we are going to use the equivalent form of the condition $(C)$, i.e., equation \eqref{b1}. Multiplying it by $\theta$ we get the equation
\beq\label{cc}
\p_x T(x,z) - T(x,z) \p_x\ln \theta (Ux+z)=0,  \quad z\in \widehat Y
\eeq
where $T(x,z):=\p_y \theta(Ux+Vy+z)|_{y=0}$. Since at $x=q(0)$ the function $\p_x\ln \theta (Ux+z)$
has simple pole with residue $N$, equation \eqref{cc} implies that $T$ at $x=q(0)$ has zero of the same order as $\theta(Ux+z)$. Therefore, the function
$\p_V\ln \theta |_{\Theta\cap\widehat Y}$ is a regular function. The latter implies that $\p_y \psi(x,y)|_{y=0}$ for $\psi$
defined in \eqref{u0} has pole $x=q(0)$ of order at most $n$. Two other terms on the left hand side of \eqref{lax0} have poles of order $n+2$. The comparison of their leading coefficients gives the equation
\beq\label{N2}
2N=n(n+1)
\eeq
From \eqref{N1} and \eqref{N2} it follows that $N=n=1$. The implication $(C)\Rightarrow (B)$ is proved. That completes the proof of the "only if" part of the theorem statement.

\medskip
We begin the proof of the "if" part of the Theorem statement by proving

\begin{lm} \label{lm:main} Suppose that equation \eqref{C} is satisfied and let $k^{-1}$ be the local coordinate defined in Lemma \ref{lm:periodicity}. Then, for the expansions \eqref{psi}, \eqref{psidual} of the Baker-Akhiezer function and its dual in this coordinate, the equation
\beq\label{duality}
\psi^*(x,0,k)=\psi(x,0,-k)
\eeq
holds.
\end{lm}
\noindent
{\it Proof.} Equation \eqref{duality} is equivalent to the equations
\beq\label{duality1}
\xi_s^*(x,0)=(-1)^s \xi_s(x,0), \ s=0,1,\ldots
\eeq
The latter, by \eqref{psiphi1} and \eqref{psiphi1star}, are corollaries of the equations
\beq\label{induc1}
\varphi_s(z)=(-1)^s \varphi_s^*(z)
\eeq
We are going to prove \eqref{induc1} by induction. The initial step of the induction is the equation \eqref{res1}.

\bigskip
Suppose that equations \eqref{induc1} hold for $s\leq n$. Then for $s\leq n$ equations \eqref{duality1} hold. The latter and equations \eqref{eqxi}, \eqref{eqxidual} with $s=n$ imply the equation
\beq\label{25}
\p_y\xi_{n}(x,0)+(-1)^n \p_y\xi_n^*(x,0)=2\p_x(\xi_{n+1}(x,0)+(-1)^{n}\xi^*_{n+1}(x,0)).
\eeq
A priory the left hand side of \eqref{25} has poles of second order at poles of $\xi_n$ and $\xi_n^*$ (i.e. at zeros of $\tau(x\eta)$) but equation \eqref{C} ensures that it has simple poles, only. At the same time the right hand side has no $x$-residues at its poles. Hence, the function $2\p_x(\xi_{n+1}(x,0)+(-1)^{n}\xi^*_{n+1}(x,0))$ is a {\it holomorphic}
function of $x$. By \eqref{psiphi1} and \eqref{psiphi1star} it equals to $f_{n+1}(x\eta)$, where
\beq\label{gf1}
f_{n+1}(z):= 2\p_\eta(\varphi_{n+1}(z)+(-1)^n\varphi^*_{n+1}(z))
\eeq
Consider the  monodromies of $f_{n+1}$. From the induction assumption and equations \eqref{phitrans1},\eqref{phitrans2} it follows  that
\beq\label{aug30}
\rho_{\l,s}=(-1)^s\rho^*_{\l,s}, \quad s=1,\ldots, n.
\eeq
The latter imply that the function $f_{n+1}$ has trivial monodromies
\beq\label{abel}
f_{n+1}(z+\lambda)=f_{n+1}(z),\quad \lambda\in \Lambda_U
\eeq
i.e. it descents to an abelian function on $Y$. Since upon restriction on the line $x\eta$ it is holomorphic, and the projection of the line is dense in $Y$ we get that $f_{n+1}(z)=f_{n+1}$ is a constant. Then the equation $\xi_s(0,0)=\xi_s^*(0,0)=0$ imply
\beq\label{feta}
F_{n+1}(x\eta)=xf_{n+1}, \quad F_{n+1}(z):=\varphi_{n+1}(z)+(-1)^n\varphi^*_{n+1}(z)
\eeq

From \eqref{aug30}, \eqref{phitrans1} and \eqref{phitrans2} it follows that
\beq\label{mon30}
F_{n+1}(z+\lambda)-F_{2n}(z)=\rho_{\lambda, n+1}+(-1)^n\rho_{\lambda,n+1}^*
\eeq
The latter imply that $F_{n}(z)$ is a linear form, i.e.
\beq \label{j12}
F_{n+1}(z)=(\tilde\ell_{n+1},z), \quad (\tilde \ell_{n+1},\eta)=f_{n+1}
\eeq
Recall, that by our choice of the local coordinate we have that $\rho_{\lambda_i, m}=\rho^*_{\lambda_i, m}=0$ for any $m>0$, i.e. the linear form $F_{n+1}$ vanishes at $d$ linear independent vectors. Hence, $F_{n+1}=0$ and the induction step is completed.

\begin{cor} Let $k^{-1}$ be the local coordinate defined in Lemma \ref{periodicity}. Then for the expansion in this coordinate  of the differential $d\Omega$ defining the dual Baker-Akhiezer function (i.e. the differential vanishing at the poles of the Baker-Akhiezer function), the equation
\beq\label{duality10}
d\Omega(k)=-d\Omega(-k)
\eeq
holds.
\end{cor}
\noindent
{\it Proof.} Equation \eqref{duality10} is equivalent to the equations
\beq\label{duality11}
c_{2s-1}=0
\eeq
where $c_s$ are the coefficients in \eqref{coeff}.

Let $c_{2n+1}\neq 0$ be the non-zero coefficient of $d\Omega$ expansion with the smallest odd index.
Let $d\Omega_{2n}$ be the sum of the first $2n$ terms of $d\Omega(k)$ expansion, i.e.,
$$d\Omega(k)=d\Omega_{2n}(k)+dk\left(c_{2n+1}k^{-2n+1}+O(k^{-n-2})\right)$$
The differential $d\Omega_{2n}$ is odd, $d\Omega_{2n}(k)=-d\Omega_{2n}(-k)$. Equation \eqref{duality} implies that $\left(\p^n_x\psi\p_x^n\psi^*_{y=0}\right)d\Omega_{2n}$ is also odd. The residue of an odd differential is zero. Then using \eqref{nomega} with $i=m=n$ we get
\beq\label{aug31}
0=\res_{k=\infty}\left(\p^n_x\psi\p_x^n\psi^*|_{y=0}\right)d\Omega=(-1)^{n+1}c_{2n+1}
\eeq
Contradiction. The corollary is proved.

\bigskip
Introduce new local coordinate $ k'=k+O(k^{-1})$ by the equation
\beq\label{kprime}
dk'=d\Omega(k)
\eeq
From \eqref{duality10} it follows that this change preserves the equation \eqref{duality}, i.e.
\beq\label{duality31}
\psi^*(x,0,k')=\psi(x,0,-k')
\eeq
For further use note that in the new coordinate the residue equation for the BA function and its dual has the standard form
\beq\label{resnew}
\res_{k=\infty}\psi^*(x,0,k')(\p_x^i \psi(x,0,k')) dk'=0
\eeq
\begin{cor}\label{cor}
Let $\LL$ be a unique pseudo-differential operator of the form \eqref{satokp}
such that the equation
\beq\label{L1}
\LL\psi(x,0,k')=k'\psi(x,0,k')
\eeq
holds. Then the equation \eqref{ckp01}, i.e.  $\LL^*=-\LL$, holds.
\end{cor}
\noindent
{\it Proof.} The proof is standard in the KP theory. Introduce, the so-called wave pseudo-differential operator $\Phi=1+O(\p^{-1})$ by the equation
\beq\label{W}
\psi(x,0,k')= \Phi e^{k'x}
\eeq
Then
\beq\label{L}
\LL=\Phi\cdot \p_x \cdot \Phi^{-1}
\eeq
It is known that \eqref{resnew} implies
\beq\label{Phidual}
\psi^*(x,0,k')=(\Phi^{-1})^* e^{-k'x}
\eeq
(see in \cite{kr-schot}). Using \eqref{duality31} we get $\Phi^*=\Phi^{-1}$. The latter and \eqref{L}) implies \eqref{ckp01}. The corollary is proved.

\medskip
Now we are going to use few more standard facts from the KP theory and the related theory of commuting ordinary differential operators.
Let $B_n=\LL^n_+$ be, as before, the differential part of the operator $\LL^n$.
From \eqref{ckp01} it follows that
\beq\label{antiself1}
B_n^*=(-1)^n B_n
\eeq
The commuting flows of the KP hierachy \eqref{kp3}
correspond to linear flows on the Jacobian. Hence, among the two sets of flows corresponding to even or odd $n$ there are only finite number of linear independent. Hence,  for all but a finite number of integers $n$ there exists a linear operator
\beq
L_n=B_n+\sum_{m=1}^{[n/2]} c_{n,m} B_{n-2m}
\eeq such that
\beq\label{aeigen}
L_n\psi(x,0,p)=a_n(p)\psi(x,0,p)
\eeq
where $a_n(p)$ is a meromorphic function on $\G$ with the only pole at the marked point  $P$.
The operator $L_n$ is self-adjoint for even $n$ and anti self-adjoint for odd $n$. In both cases we
have
\beq\label{aeigen1}
L_n\psi^*(x,0,p)=a_n^*(p)\psi^*(x,0,p)
\eeq
where $a_n^*(p)$ is a meromorphic function of $\G$ with the only pole at $P$ of order $n$. Since  $\psi^*(x,0,k')=\psi(x,0,-k')$, the expansions of $a_n$ and $a_n^*$ near $P$ in $k'$ are related by the equation
$$a_n^*(k')=(-1)^n a_n(-k')
$$
The involution $a_n\to a^*_n$ of the ring $\mathcal A(\G,P)$ of meromorphic functions on $\G$ with the only pole at $P$ defines the involution $\s:\G\to \G$ under which the local parameter $k^{-1}(p)$ is odd, $k'(\sigma (p))=-k'(p)$. The theorem is proved.

\section{Proof of Theorem \ref{main2}} \label{s:2d}

To begin with, recall necessary facts form the algebraic-geometrical construction of solutions of the $2D$ Toda hierarchy.

\medskip
\subsection*{Two-point Baker-Akhiezer function}

Let $\G$ be a smooth genus $g$ algebraic curve with fixed coordinates $k_\a^{-1}(p)$ in the neighborhoods of two marked points $P_\a,\,\a=1,2$. The two-point Baker-Akhiezer function $\psi$ depends on the variables $(x, \bt=(t_{\a,i>0}))$ (as always, it is assumed that only a finite number of them are not zero). For fixed $(x,\bt)$ and a non-special effective degree $g$ divisor $D$ it is defined as a function of $p\in \G$ with the following analytical properties.

On $\G\setminus \{P_1,P_2\}$ the function $\psi$ is a multi-valued meromorphic function with poles at $D$ whose order at a point of $D$ is not greater than the multiplicity of the point in $D$, i.e. $(\psi)+D\geq 0$. The multi-valuedness of $\psi$ is due to the only not trivial monodromies around the marked points $P_1$ and $P_2$ which are equal to $e^{2\pi ix}$ and $e^{-2\pi i x}$, respectively (for integer $x$ the function $\psi$ is single-valued).

In the neighborhood of $P_1$ it has the form:
\beq\label{psiP1}
\psi (x,\bt ,p) =
k_1^x\exp \biggl(\sum_{i=1}^{\infty} t_{1,i} k_{1}^{i} \biggr)
\biggl(1+ \sum_{s=1}^{\infty} \xi_{1,s}(x,\bt ) k_1^{-s}\biggr),
\eeq
and in the neighborhood of $P_2$
\beq\label{psiP2}
\psi (x,\bt ,p) =
k_2^{-x}e^{f(x,\bt)}\exp \biggl(\sum_{i=1}^{\infty} t_{2,i} k_{2}^{i} \biggr)
\biggl(1+\sum_{s=1}^{\infty} \xi_{2,s}(x,\bt) k_1^{-s}\biggr),
\eeq
As shown in \cite{kr-toda} the Baker-Akhiezer function satisfies the linear equations
\beq\label{nov23}
(\p_{1,1}-T-u)\psi=0,\ \ \ (\p_{2,1}-w T^{-1})\psi=0,
\eeq
where
\beq\label{nov231}
u=\p_{1,1} f=(1-T)\, \xi_{1,1}, \ \ \ w=\exp((1-T^{-1})f),\ \
\eeq
The compatibility condition \eqref{nov23} is equivalent to the $2D$ Toda equation
\beq\label{2DT}
\p_\xi\p_\eta f_n=e^{f_{n}-f_{n-1}}-e^{f_{n+1}-f_{n}},
\eeq with $x=n, \xi=t_{1,1}$ and $\eta=t_{2,1}$.

\medskip
The explicit formula for $f$ follows from the
theta-functional formula for the \BA functions:
\beq\label{BA2} 
\psi=
\frac{\theta(A(P)+xU+\sum U_{\a,i}t_{\a,i}+Z) \,\theta(Z)}
{\theta(xU+\sum U_{\a,i}t_{\a,i}+Z)\, (\theta(A(P)+Z) }
e^{\left(x\Omega_0+\sum t_{\a,i}\Omega_{\a,i}(P)\right)},
\eeq
Here $\Omega_0$ is an Abelian integral of the normalized meromorphic differential with simple poles at $P_1,P_2$ and residues $-1$ and $1$ respectively. From the bilinear Riemann identities it follows that
\beq\label{U2d}
U=A(P_2)-A(P_1)
\eeq
where $2\pi i U$ is the vector of $b$-periods of $d\Omega_0$. The Abelian integrals $\Omega_{\a,i}$ and the corresponding vectors $U_{\a,i}$ are defined similarly to that in \eqref{formula}: $d\Omega_{\a,i}$ is the normalized meromorphic differential with the pole at $P_\a$ of the form $d(k^i_{\a}+O(k_\a^{-1}))$. The corresponding Abelian integrals are normalized such that their expansions near $P_1$ are of the form:
\beq\label{const1}
\Omega_{1,i}=k^i+\sum_{s=1}^\infty a_{1,s}^{(1,i)} k_1^{-s}\,,  \quad \Omega_{2,i}=\sum_{s=1}^\infty a_{1\,s}^{(2,i)}k_1^{-s}
\eeq
and near $P_2$
\beq\label{const2}
\Omega_{1,i}=\sum_{s=0}^\infty a_{2,s}^{(1,i)} k_2^{-s}\,,  \quad \Omega_{2,i}=k_2^i\sum_{s=1}^\infty a_{2\,s}^{(2,i)}k_2^{-s}
\eeq
The expansions of $\Omega_0$ are of the form:
\beq\label{caonst0}
\Omega_0=\ln k_1+\sum_{s=1}^\infty a_{1,s}^{(0)} k_1^{-s}, \quad \Omega_0=-\ln k_2+\sum_{s=1}^\infty a_{2,s}^{(0)} k_2^{-s}
\eeq
In the definition of the Abel transform we chose the normalization $A(P_1)=0$. Then using \eqref{U2d} we get
\beq\label{ftau}
f=(T-1)\ln \theta +x a_{2,0}^{(0)}+\sum t_{\a,i} a_{2,0}^{(\a,i)},
\eeq
$$\theta:=\theta(xU+\sum U_{\a,i}t_{\a,i}+Z)$$

For further use notice that \eqref{BA2} and the expansion \eqref{Ak} of the Abel
transform near $P_1$ imply
\beq\label{xi11}
\xi_{1,1}=-\p_{1,1} \ln\theta +x a_{1,1}^{(0)}+\sum t_{\a,i} a_{1,1}^{(\a,i)}
\eeq

\medskip

In what follows we restrict ourselves to the case the first three variables $x,t_{1,1},t_{2,1}$ by setting
$t_{\a,i>1}=0$. Introduce the variables $y$ and $t$ such that $t_{11}=y+t,\, t_{2,1}=-t$.

In the new variables equations \eqref{nov23} take the form:

\begin{lm} The Baker-Akhiezer function $\psi(x,y,t,p)$ given by formula \eqref{BA2} satisfies the equation
\beq\label{laxy}
\left(\p_y-T-u\right)\psi=0,\ \ T=e^{\p_x}
\eeq
with
\beq\label{ud}
u=b_1+(T-1)\p_y \ln \theta, \quad b_1:=a_{1,1}^{(0)}
\eeq
and the equation
\beq\label{laxt}
(\p_t-T-u-wT^{-1})\psi=0
\eeq
with
\beq\label{vwd}
w=b_2\frac{T\theta \, T^{-1}\theta}{\theta^2}, \quad b_2:=e^{a_{2,0}^{(0)}}
\eeq
\end{lm}
The substitution of \eqref{psiP1} into \eqref{laxy} gives a system of the equations
\beq\label{eqxi11}
\p_y \xi_{1,s}=(T-1)\xi_{1,s+1}+u\xi_{1,s}, \quad s=0,1,2,\ldots
\eeq
The first of them gives an expression of $u$ in terms of $\xi_{1,1}$
\beq\label{xi1u}
(1-T)\,\xi_{1,1}=u=b_1+(T-1)\p_y\ln \theta
\eeq
Similarly, the substitution of \eqref{psiP2} into \eqref{laxy} gives:
\beq\label{eqxi2}
\p_y \xi_{2,s}=e^{(T-1)f} \left(T\xi_{2,s-1}\right)=b_2\,\frac{T^2\theta \cdot \theta}{(T\theta)^2}\,\left(T\xi_{2,s-1}\right), \quad s=0,1,2, \ldots
\eeq
where for the last equation we use \eqref{ftau}.

\medskip

For further use introduce the function
\beq\label{dphi}
\phi(z, k_2):= \frac{\theta({A}(k_2)+\iota(z))\,\theta(Z)}
{\theta(\iota(z))\, \theta(A(k_2)+Z)}\, e^{(\ell(k),z)}\,, \quad z\in \mathbb C^d,
\eeq
where $A(k_2)$ is the expansion of the Abel transform near $P_2$ and $\ell(k_2)$ is the formal series
\beq\label{dell}
\ell(k_2)=\sum_{s=0}^\infty \ell_{s}k_2^{-s}, \quad \ell_{s}\in \mathbb C^d,
\eeq
such that the equation
\beq\label{doell}
(\ell(k_2),\eta)=\Omega_0(k_2)+\ln k_2\,,
\eeq
where $\Omega(k_2)$ is the expansion of the Abelian integral $\Omega_0$ at $P_2$, holds.

\smallskip
A simple comparison of \eqref{BA2} and \eqref{dphi} shows that the equation
\beq\label{dpsiphi1}
\xi_s(x,0)=\varphi_s(x\eta)
\eeq
where $\varphi_s$ are the coefficients of the expansion of $\phi$.
The latter are of the form
\beq\label{phi111}
\varphi_s(z)=\frac {\tau_{2,s}(z)}{\tau(z+\eta)} , \quad \tau(z):=\theta(\iota(z))
\eeq
where $\tau_{2,s}(z)$ is a holomorphic function of $z$.

The monodromy properties of the theta-function imply that under the shift of the arguments by a vector $\lambda\in \Lambda_U$  the function $\phi$ gets transformed to
\beq\label{dphitrans}
\phi(z+\lambda,k_2)=\phi(z,k_2)\rho_\lambda(k_2), \ \ \quad \rho_\lambda(k_2)=1+\sum_{s=1}^\infty\rho_{\lambda,s}k_2^{-s}
\eeq
Equation \eqref{phitrans} is equivalent to the equations
\beq\label{dphitrans1}
\varphi_s(z+\lambda)=\sum_{j=0}^s \varphi_{s-j}(z)\rho_{\lambda, j}
\eeq
The same arguments as in the proof of Lemma \ref{lm:periodicity} prove:
\begin{lm}\label{lm:dperiodicity} Let $\l_1,\ldots,\l_d$ be a set of linear independent vectors in $\Lambda_U$. Then there is a unique coordinate $k_2^{-1}(p)$ in the neighborhood of $P_2$ and a unique  linear form $\ell(k_2)$ satisfying \eqref{doell} such that $\phi(z,k_2)$ is invariant under the shift by these vectors, i.e.,
\beq\label{periodicity2}
\phi(z+\lambda_i,k_2)=\phi(z,k_2) \quad \Leftrightarrow \quad \rho_{\lambda_i}(k_2)=1, \ \ i=1,\ldots, d.
\eeq
\end{lm}
\subsection*{Dual Baker-Akhiezer function}

As in the one-point case the definition of the dual two-point Baker-Akhiezer function begins with the definition of the dual divisor. Let $D$ be an effective degree $g$ divisor. In the two point case the dual divisor $D^*$ is defined as the divisor such that $D+D^*$ is the zero divisor of the meromorphic differential $d\Omega$ with simple poles at the marked points, i.e.,
\beq\label{Dstar}
D+D^*=K+P_1+P_2\in J(\G)
\eeq
In what follows it will assumed that $d\Omega$ is normalized such that it has residue $-1$ at $P_1$.

The dual Baker-Akhiezer function $\psi^*(x,\bt,p)$ as a function of $p\in \G$ has the following analytic properites: away of the marked points it is multi-valued meromorphic function with the pole divisor $D^*$. Its the only non-trivial monodromy factors are $e^{-2\pi ix}$ around $P_1$ and $e^{2\pi i x}$ around $P_2$.

In the neighborhood of $P_1$ it has the form:
\beq\label{dpsiP1}
\psi^* (x,\bt,p) =
k_1^{-x}\exp \biggl(\sum_{i=1}^{\infty} -t_{1,i} k_{1}^{i} \biggr)
\biggl(1+ \sum_{s=1}^{\infty} \xi^*_{1,s}(x,\bt) k_1^{-s}\biggr),
\eeq
and in the neighborhood of $P_2$
\beq\label{dpsiP2}
\psi^* (x,\bt,p) =
k_{\,2}^{x}\,e^{f^*(x,\bt)}\exp \biggl(\sum_{i=1}^{\infty} -t_{2,i} k_{2}^{i} \biggr)
\biggl(1+\sum_{s=1}^{\infty} \xi^*_{2,s}(x,\bt) k_1^{-s}\biggr),
\eeq
The explicit theta-functional formla for th
\beq\label{BA2s} 
\psi^*=\frac{\theta(A(P)-(x+1)U-\sum U_{\a,i}t_{\a,i}-Z) \,\theta(Z-U)}
{\theta((x+1)U+\sum U_{\a,i}t_{\a,i}+Z)\, (\theta(A(P)-Z-U) }e^{-x\Omega_0-\sum t_{\a,i}\Omega_{\a,i}(P)}
\eeq
where we use the equality $Z^*=-Z-U$ which follows from \eqref{Dstar} and \eqref{U2d}.

From the definition of the dual divisor it is easy to see that the equations
\beq\label{resd}
\res_{P_1}\left(\psi^* T^i \psi\right) d\Omega= -\delta_{i,0}
\eeq
hold. Indeed  the differential in the r.h.s. of \eqref{resd} is a meromorphic differental on $\G$ which is holomorphic away of the marked points. For $i>0$ it is holomorphic at $P_2$. Hence, has no residue at $P_1$. Equation \eqref{resd} for $i=0$ follows from \eqref{psiP1} and \eqref{dpsiP1} and the normalisation of $d\Omega$ defined above.

Notice, that \eqref{resd} with $i=0$ implies the equation
\beq\label{resd0}
f(x,\bt)=-f^*(x,\bt)
\eeq
and with $i=1$ implies the equation
\beq\label{indd1}
\xi_{1,1}(x+1,\bt)+\xi_{1,1}^*(x,\bt)=0
\eeq
Equation \eqref{resd0} implies that the dual BA function satisfies the equation formally adjoint to \eqref{laxy}
\beq\label{laxys}
(\p_y+T^{-1}+u)\psi^*=0
\eeq
The substitution of \eqref{dpsiP1} into \eqref{laxys} gives the equations
\beq\label{eqxi1d}
-\p_y\xi_{1,s}^*=(T^{-1}-1)\xi^*_{1,s+1}+u \xi_{1,s}^*, \quad s=0,1,\ldots
\eeq

\medskip
For further use introduce the function
\beq\label{dphist}
\phi^*(z, k_1):= \frac{\theta({A}(k_1)-\iota(z+\eta))\,\theta(Z+U)}
{\theta(\iota(z+\eta))\, \theta(A(k_1)-Z-U)}\, e^{-(\ell_1(k_1),z)}\,, \quad z\in \mathbb C^d,
\eeq
where $A(k_1)$ is the expansion of the Abel transform near $P_1$ and $\ell_1(k_1)$ is the formal series
\beq\label{dellst}
\ell(k_1)=\sum_{s=1}^\infty \ell_{1,s}k_1^{-s}, \quad \ell_{1,s}\in \mathbb C^d,
\eeq
such that the equation
\beq\label{doellst}
(\ell_1(k_1),\eta)=\Omega_0(k_1)-\ln k_1\,,
\eeq
where $\Omega_0(k_1)$ is the expansion of the Abelian integral $\Omega_0$ at $P_1$, holds.

A simple comparison of \eqref{BA2s} and \eqref{dphist} shows that the coefficients of the dual Baker-Akhiezer function
expansion at $P_1$ are equal to
\beq\label{dpsiphi1s}
\xi^*_{1,s}(x,0)=\varphi^*_s(x\eta)
\eeq
where $\varphi^*_s$ are the coefficients of the expansion of $\phi^*$.
The latter are of the form
\beq\label{phi11d}
\varphi^*_s(z)=\frac {\tau^*_{1,s}(z)}{\tau(z+\eta)} , \quad \tau(z):=\theta(\iota(z))
\eeq
where $\tau^*_{2,s}(z)$ is a holomorphic function of $z$.

The monodromy properties of the theta-function imply that under the shift of the arguments by a vector $\lambda\in \Lambda_U$  the function $\phi^*$ gets transformed to
\beq\label{dphitransst}
\phi^*(z+\lambda,k_1)=\phi(z,k_2)\rho^*_\lambda(k_1), \ \ \quad \rho^*_\lambda(k_1)=1+\sum_{s=1}^\infty\rho^*_{\lambda,s}k_1^{-s}
\eeq
Equation \eqref{dphitransst} is equivalent to the equations
\beq\label{dphitrans1st}
\varphi^*_s(z+\lambda)=\sum_{j=0}^s \varphi^*_{s-j}(z)\rho^*_{\lambda, j}
\eeq
The same arguments as in the proof of Lemma \ref{lm:periodicity} prove:
\begin{lm}\label{lm:dperiodicitystst} Let $\l_1,\ldots,\l_d$ be a set of linear independent vectors in $\Lambda_U$. Then there is a unique coordinate $k_1^{-1}(p)$ in the neighborhood of $P_1$ and a unique  linear form $\ell_1(k_1)$ satisfying \eqref{doellst} such that $\phi^*$ is invariant under the shift by these vectors, i.e.,
\beq\label{dstperiodicity}
\phi^*(z+\lambda_i,k_1)=\phi^*(z,k_1) \quad \Leftrightarrow \quad \rho^*_{\lambda_i}(k_1)=1, \ \ i=1,\ldots, d.
\eeq
\end{lm}

\bigskip
\subsection*{ Proof of the theorem.} We begin the proof of the theorem by pointing on the obstacle which makes our approach for the characterization of the Jacobian of curves with involution non-applicable for the case of non-ramified covers.

\begin{lm} \label{lm:1}Let $\G$ be a curve with involution having fixed points. Then for the differential $d\Omega_0$ defined above, i.e., the Abelian differential with simple poles at the marked points $P_1,P_2$, which are permuted by involution, and normalized by the conditions $\oint_{a_k} d\Omega_0=0$, the equation
\beq\label{difodd}
d\Omega_0(\s(p))=-d\Omega_0(p),\quad p\in \G
\eeq
holds.
\end{lm}
\begin{rem}{\rm
The author have not been able to find a  reference to this statement in the literature, but he has no doubt that it is well-known. In any case we present a proof of it, since its arguments will be useful below.}
\end{rem}

\noindent
{\it Proof.} Following the line of the proof of Lemma 3.3 in \cite{kr-quad}  it is easy to show that there are integers $l_k$ such that the equation
\beq\label{difodd1}
d\Omega_0(\s(p))=-d\Omega_0(p)+2\pi i\, l_k \,\omega_k\,, \quad l_k\in \mathbb Z\,,
\eeq
where $\omega_k$ are normalized holomorphic differentials, holds.

On any compact  $\bar\G\in\{\G\setminus {P_1,P_2}\}$ the differential $d\Omega_0$ depends continuously on  the marked points $P_1,P_2=\s(P_1)$. Hence, integers $l_k$ do not depend on $P_1$. In the limit $P_1\to P$, where $P$ is a fixed point of involution, the residues of $d\Omega_0$ cancel each other. The latter implies that on $\bar\G$ the differential $d\Omega_0$ tends to zero uniformly.  Then comparing the limits of two sides of \eqref{difodd1} we get  $l_k=0$. The lemma is proved.

\begin{rem} {\it \rm The statement of the lemma above fails to be true for the case of unramified covers.}
\end{rem}

\begin{lm} \label{lm:2}Suppose that there is an involution $\s:\G\to \G$ with fixed points under which the marked points are permuted, $\s(P_1)=P_2$ and
the local coordinates in the neighborhoods of marked points are chosen such that $k_1(p)=k_2(\s(p))$.  Then the equation
\beq\label{turning}
\left(\theta^2\p_y^2\ln\theta-b_2T\theta \,T^{-1}\theta \right)\Big|_{y=0}=b_3
\eeq
where $b_3$ is a constant and $\theta=\theta(xU+yV+tW+\zeta)$ with $\zeta$ such that $\zeta+\zeta^\s=K+P_1+P_2$, where $K$ is the canonical class,
holds.
\end{lm}
\noindent
{\it Proof.} Let $\psi$ be the BA function defined by a divisor $D$ such that
$D+D^\s$ is the zero divisor of a meromorphic diffrential $d\wt\Omega$ with simple poles at $P_1$ and $P_2$. The vector $\zeta$ above is equal to $\zeta=-A(D)+\mathcal K$, where $\mathcal K$ is the vector of Riemann constants.

\medskip
Consider the differential $d\wh\Omega:=(T\psi)\psi^\s|_{y=0}d\wt\Omega$. It is holomorphic outside the marked points, since poles of $\psi$ and $T\psi^\s$ cancel out with zeros of $d\wt\Omega$. At the marked points essential singularities of $T\psi$ and $\psi^\s$ cancel each other. The multi-valuedness of the BA function requiers additional consideration.
Notice  that the differential $d\Omega_{1,1}-d\Omega_{2,1}$ is odd with respect to the involution. Then from \eqref{BA2} and \eqref{difodd} it follows that $d\wh\Omega$ is {\it single-valued}, i.e. it is a meromorphic differential on $\G$ with the only pole at $P_1$. Hence, it has no residue at that point. Computing the residue  in terms of the coefficients of expansions \eqref{psiP1} and \eqref{psiP2} we get the equation
\beq\label{a26}
\xi_{1,1}(x+1,0,t)+\xi_{2,1}(x,0,t)=0
\eeq
Then from equation \eqref{eqxi2} with $s=1$ and equation \eqref{xi11} it follows that equation \eqref{turning} holds. The lemma is proved.

Equation \eqref{turning} implies equation\eqref{bd1}, since by definition $\pi(Ux+Wt+\zeta)$ is dense in $\wh Y$.
Note also, that equation \eqref{turning} is equivalent to equation \eqref{Cd}.
Indeed, substituting \eqref{psiP2} into \eqref{laxt} and comparing the leading coefficients of the expansions at $P_2$ we get
\beq\label{utt}
u=\p_t f+(T^{-1}-1)\xi_{2,1}=\p_t f+(T-1)\xi_{1,1}=\p_t f-u
\eeq
where  in the last two equalities we use \eqref{a26} and \eqref{nov231}, respectively.
From \eqref{utt} it follows that
\beq\label{utt1}
2u=2\p_y f|_{y=0}=\p_t f
\eeq
Equation \eqref{Cd} is proved.

\medskip
Our next goal is to prove the condition $(B)$ of the theorem. The explicit meaning of part $(i)$ of $(B)$ is as follows: generic zero $x=q$ of $\tau(x\eta+z)$ with $z\in \pi^{-1}(Y)$ is simple,  and
$\tau((q+1)\eta+z)\tau((q-1)\eta+z)\neq 0$.
Here and below $\tau(z,y,t):=\theta(\iota(z)+Vy+Wt)$, and for brevity we will write $\tau(z,y):=\tau(z,y,0)$ and
$\tau(z,t):=\tau(x,0,t)$, and $\tau(z)=\tau(z,0,0)$.
\begin{rem} {\rm As it was mentioned in Introduction the author was unable to complete a proof that equations \eqref{Cd} and \eqref{bd1} imply $(B)$ although belives that it is true. What we were not able to exclude is the possibility that $\Theta\cap Y$ is a union of two components $(\Theta\cap Y)_\pm$ on which zeros of $\tau$ are of multiplicity $2$, and $(\Theta\cap Y)_-+U=(\Theta\cap Y)_+$. The proof below uses an extension of arguments in the proof of Lemma \ref{lm:1} and the result of the previous Section.}
\end{rem}

Consider in more details the differential $d\Omega_0$ when $P_1$ is close to a fixed point $P$ of the involution.
Let $z$ be an odd with respect to $\s$ coordinate in the neighborhood of $P, \, z(P)=0$. Let $z_0=z(P_1)$ then in the neighborhood of $P$ the differential $d\Omega_0$ has the form
\beq\label{lim1}
d\Omega_0=\frac {dz}{z+z_0}-\frac {dz}{z-z_0}+R(z,z_0)dz
\eeq
where $R(z,z_0)$ is even holomorphic function of $z$ (in the neighborhood of $z=0$) and odd function of $z_0$ (uniformly bounded by $z_0$). Equation \eqref{lim1} obviously shows that $d\Omega_0\to 0$ as $z_0\to 0$. More importantly, it shows that $(2z_0)^{-1} d\Omega_0\to d\Omega_1$, where $d\Omega_1$ is the normalized differential with the pole of order $2$ at $P$ which we used in the previous section.

\medskip\noindent
{\it Change of notations.} For the duration of these arguments we use notation $2\pi i U_0$ and $2\pi i U_1$ for the vectors of $b$-periods of the differentials $d\Omega_0$ and $d\Omega_1$ in order to avoid a confusion caused by use the same notation $U$ for them in Theorems \ref{thm:main} and \ref{main2}. Similarly, we will use the notations
$Y_0$ and $Y_1$ for the closers of $\pi(U_0x+\zeta_0)$ and $\pi(U_1x+\zeta_1)$ in $J(\G)$, respectively (which both  were denoted above by $Y$).

\medskip
The convergence $(2z_0)^{-1} d\Omega_0\to d\Omega_1$ is uniform on any compact $\bar \G\subset (\G\setminus P)$. The latter implies the convergence of the vectors $(2z_0)^{-1} U_0\to U_1$ and as a corollary the convergence
$Y_0\to Y_1$.

As it was shown in the proof of Theorem \ref{thm:main} the intersection $\Theta\cap Y_1$ is reduced. The latter condition is open. Hence, for $P_1$ sufficiently closed to $P$ the $\Theta\cap Y_0$ is also reduced. Then, using again openness of the latter condition, we conclude that $\Theta\cap Y_0$ is reduced for all but at most  finite number of points $P_1\in \G$. Upper continuity of multiplicty of zeros of entire function depending on parameters implies that $\Theta\cap \wh Y_0$, where $\wh Y_0$ is the closure of $\pi(U_0x+Wt+\zeta_0)$, is also reduced.

If $\tau(x,0,t)$ has simple zero at $x=q(t)$, then from equations \eqref{turning} and \eqref{utt1} it follows
the equations
\beq\label{qdot}
\left(\p_t q(t)\right)^2=-2b_2\,\frac{\tau(q+1,0,t)\tau(q-1,0,t)}{\tau^2_x(q,0,t)}
\eeq
Without loss of generality (see Remark \ref{rm:W}) we may assume that the vector $W\neq 0$. Hence, generically the left hand side of \eqref{qdot} is not zero. Hence, generically, the r.h.s.of the equation does not vanish. That is the statement $(B)$ of the theorem. The "only if" part of the theorem is proved.

\bigskip
The following lemma is central in the proof of the "if" part of the theorem.

\begin{lm} \label{lm:maind} Suppose that equation \eqref{Cd} is satisfied. Then for the coefficients of the expansions\eqref{psiP2}, \eqref{dpsiP1} of the Baker-Akhiezer function and its dual in the coordinates $k_1^{-1}$ and $k^{-1}_2$ in the neighborhoods of the marked points defined in Lemmas \ref{lm:dperiodicity} and \ref{lm:dperiodicitystst} for the same set of $d$ independent vectors $\lambda_i\in \Lambda_U$ the equations
\beq\label{dduality1}
\xi_{1,s}^*(x,0)=\xi_{2,s}(x,0), \ s=0,1,\ldots
\eeq
hold.
\end{lm}
\noindent
{\it Proof.} By \eqref{dpsiphi1} and \eqref{dpsiphi1s} equations \eqref{dduality1} are corollaries of the equations
\beq\label{induc1d}
\varphi_s(z)=\varphi_s^*(z), \quad s=0,1, \ldots
\eeq
We are going to prove \eqref{induc1d} by induction. The initial step of the induction is the equation $\varphi_0=\varphi_0^*=1$. Suppose that equation \eqref{induc1d} is satisfied for $s\leq n$.

Compute the residue $r_{n+1}$ of $\xi_{2,n+1}(x,0)$  at the zero $q$ of $\tau(x+1)$. Comparing the coefficient at $(x-q)^{-2}$ of the left and right sides of equation \eqref{eqxi2}
we get
\beq\label{rn}
-((\p_y q)|_{y=0}) r_{n+1}=b_2\frac{\tau(q+2)\tau(q)}{\tau^2_x(q+1)}\xi_{2,n}(q+1,0)
\eeq
Compute now the residue $r^*_{n+1}$ of $\xi^*_{1,n+1}$. Consider the equation \eqref{eqxi1d}. The l.h.s. of the equation and the second term of the r.h.s. are regular at $x=q+1$. Two remaining terms have simple poles at $x=q+1$. Comparing their residues at $x=q+1$ we get the equation
\beq\label{rns}
r_{1,n+1}^*=((\p_y q)|_{y=0})\xi^*_{1,n}(q+1,0)
\eeq
By induction assumption $\xi_{2,n}(q+1,0)=\xi^*_{1,n}(q+1,0)$. Therefore,  equations \eqref{rn}, \eqref{rns} and \eqref{qdot} imply
\beq\label{rnf}
r_{2,n+1}=r^*_{1,n+1}
\eeq
From \eqref{rnf} it follows that the function $F_{n+1}(z):=\varphi_{n+1}(z)-\varphi^*_{n+1}(z)$ restricted to the line $z=x\eta$ is holomorphic. Then the same arguments as in the previous section prove that $F_{n+1}=0$. The induction step is completed. The lemma is proved.

\medskip
Note that if we identify the local coordinates $k_1^{-1}$ and $k_2^{-1}$ then equations \eqref{dduality1} take the form
\beq\label{dduality}
\psi_2(x,0,k)=e^{f(x,0)}\psi_1^*(x,0,k)
\eeq

The final steps in the proof of Theorem \ref{main2} are literally identical to those in the proof of Theorem
\ref{thm:main}, after changing the facts from KP theory on their analogues in $2D$ Toda theory. Let us briefly outline them.

First, introduce the pseudo-difference operators
\beq\label{LL1}
\LL_1=\sum_{i=-1}^\infty u_{1,i} T^{-i} ,  \quad \LL_2=\sum_{i=-1}^\infty u_{2,i} T^{i}
\eeq such that the equations
\beq\label{LLpsi}
\LL_1\, \left(e^{-f/2}\psi(x,0,k_1)\right)=k_1e^{-f/2}\psi(x,0,k_1),
\eeq
$$\LL_2\left(e^{-f/2}\psi(x,0,k_2)\right)=k_2 e^{-f/2}\psi(x,0,k_2)$$
hold. Then from \eqref{dduality1} or equivalently \eqref{dduality} it follows that the equation
\beq\label{L1L2}
\LL_1=\LL_2^*
\eeq
holds. Recall, that in the difference case the formal adjoint operator is defined by the rule $(w\cdot T^i)^*=T^{-i}\cdot w$ and then extended by linearity.

Define difference operators
\beq\label{sep16}
B_{1,m}=(({\LL_1}^m)_{>0}+\frac{1}{2}\, ((\LL_1)^m)_{0} , \quad
B_{2,m}=(\LL_2)^m)_{<0}+\frac{1}{2}\, ((\LL_2)^m)_{0}.
\eeq
where $(\cdot)_{> 0}, (\cdot)_{<0} $ denote the positive and negative parts of the pseudodifferential operator, respectively, and $ \cdot)_{0} $ denotes its zero-order term.

\medskip
The equations of $2D $ Toda hierarchy (in the so-called symmetric gauge, see \cite {kz2} for details)
\beq\label{flows2d}
\p_{\a,m} \LL_\b=[B_{\a,m} ,\LL_\b]\,, \quad \a, \b=1,2, \ \ m=1,2\ldots
\eeq
define commuting flows on the space of pairs of pseudo-difference operators of the form \eqref{LL1}.
These flows
correspond to linear flows on the Jacobian. Hence, among them there is only finite number of linear independent ones. Therefore,  for all but a finite number of integers $n$ there are constants $c_{n,m}^\pm$ such that the corresponding linear combination of these flows defined by the operators
\beq\label{LLpm}
L_n^{\pm}=\sum_{i=1}^m c_{n,m}^\pm \left(B_{1,m}\,\pm \,B_{2,m}\right)
\eeq
is trivial. From \eqref{L1L2} it follows that $(L_n^\pm)^*=\pm L_n^\pm$. The latter implies the equations
\beq\label{aeigend}
L_n^{\pm}\psi(x,0,p)=a_n^{\pm}(p)\psi(x,0,p)\,, \quad L_n^{\pm}\psi^*(x,0,p)=\pm a_n^\pm (p)\psi^*(x,0,p)\,,
\eeq
where $a_n^\pm(p)$ are meromorphic functions on $\G$ with the only poles of order $n$ at the marked points  $P_1, P_2$, hold. Points of the curve $\G$ parameterise common eigenfunctions of the commuting difference operators $L_n^{\pm}$ (see \cite{mum,kr-dif}). Hence, the correspondence $\psi\to \psi^*$ defines an involution $\s$ of $\G$ under which the equation
\beq\label{final}
\psi(x,0,p)=\psi^*(x,0,\s(p))
\eeq
hold. The theorem is proved.

\end{document}